\documentclass[12pt,a4paper]{article}
 
\usepackage[latin1]{inputenc}
\usepackage[english]{babel}
\usepackage[T1]{fontenc}
\usepackage{amsmath} 
\usepackage{amsfonts}
\usepackage{amssymb}
\usepackage{amsthm}
\usepackage{ dsfont }
\usepackage{ stmaryrd }
\usepackage{tikz}
\usetikzlibrary{matrix,arrows.meta}
\usepackage[title]{appendix}

\allowdisplaybreaks

\usepackage{hyperref}

\usepackage{xcolor} 
\usepackage[a4paper,top=2.1cm,bottom=2.60cm,left=2.1cm,right=2.1cm]{geometry} 

\theoremstyle{plain}
\newtheorem{theorem}{Theorem}[section]
\newtheorem{lemma}[theorem]{Lemma}
\newtheorem{proposition}[theorem]{Proposition}

\theoremstyle{definition}

\begin{document}

\author{Paolo Luzzini \thanks{ EPFL, SB MATH SCI-SB-JS, Station 8, CH-1015 Lausanne, Switzerland. E-mail:
paolo.luzzini@epfl.ch} }

\title{A mapping property of the heat volume potential}
\maketitle

\abstract{We consider the volume potential associated with the heat operator and we prove 
a mapping property in the space of distributions which are the time derivative of H\"older continuous functions. As an application we solve the Dirichlet and Neumann problems for the heat equation  with a non-homogeneous term in such space of distributions.}  

\vspace{10pt}

\vspace{9pt}


  \section{Introduction}
  
  One of the most useful tools to deal with non-homogeneous equations  is of course the volume potential.
  For this reason
 many authors have investigated the mapping properties of operators of volume potential type  in different functional settings and for several partial differential operators. 
 While the elliptic framework is better understood, 
 parabolic volume potentials are less investigated. For example, it
 is well know that if $F$ is a $C^{0,\alpha}$-vector field defined on a sufficiently regular bounded open subset  of 
 $\mathbb{R}^d$ then the Newtonian volume potential $\tilde P[\mathrm{div} F]$, {\it i.e.} the volume potential 
 associated with the Laplace operator and applied on $\mathrm{div}F$,  is of class $C^{1,\alpha}$ (see, {\it e.g.}, Dalla Riva, Lanza de Cristoforis and Musolino \cite{DaLaMu21}). This property allows to use the Newtonian volume potential $\tilde P[\cdot]$ to deal with the Poisson equation when the non-homogeneous term is the distributional divergence of a $C^{0,\alpha}$-vector field $F$, that is 
 \begin{equation*}\label{eq:ellF}
 \Delta u = \mathrm{div} F.
 \end{equation*}
 The parabolic analog of the above equation is 
 \begin{equation}\label{eq:heatG}
 \partial_t u-\Delta u = \mathrm{div}\, G.
 \end{equation}
Boundary value problems for equation \eqref{eq:heatG} under several H\"older regularity assumptions on the 
vector field $G$ were considered in Lieberman \cite{Li96} and Lunardi and Vespri \cite{LuVe91} with different 
techniques.
 Instead, to the best of the author knowledge, the classical parabolic theory does not cover the case
 \begin{equation}\label{eq:heatf}
 \partial_t u -\Delta u = \partial_t f,
 \end{equation}
 where $f$ is a $\frac{1+\alpha}{2}$-H\"older continuous in time and  $\beta$-H\"older continuous in space. Motivated by the above example, in the present paper we develop a theory for the volume potential associated with the heat operator acting 
 on the space of distributions of the form $\partial_t f$. As a consequence, we show how to solve the Dirichlet and Neumann  problems for equation \eqref{eq:heatf}. We note that in principle one could also try to deal with \eqref{eq:heatf} with a semigroup approach following, {\it e.g.}, Lunardi \cite{Lu95}. However, our aim is to consider 
\eqref{eq:heatf} from the point of view of potential theory and develop some tools that we plan to exploit to analyze 
perturbation problems for the heat equation via potential theory.
  
 For the classical results on elliptic volume potentials we mention here Gilbarg and Trudinger \cite{GiTr83} and Miranda 
 \cite{Mi65}. We also note that a potential theoretic approach has recently revealed 
 to be very effective to deal with elliptic problems in singularly perturbed domains. For this reason, mapping properties of elliptic volume potentials  have been also considered in view of applications to perturbation problems (see Dalla Riva,  Lanza de Cristoforis and Musolino \cite{DaLaMu15, DaLaMu20}). More details on the potential theoretic
 approach to perturbation problems for elliptic equations and results on volume potentials can be found in the monograph by Dalla Riva, Lanza de Cristoforis 
 and Musolino \cite{DaLaMu21}.

  For what concerns the parabolic case, 
 regularity properties of the heat volume potential have been considered in Friedmann \cite{Fr64}. 
Lady\v{z}enskaja,  Solonnikov, and Ural'ceva \cite{LaSoUr68} proved a series of mapping properties  of the heat volume potential in parabolic 
Schauder  and Sobolev spaces.  In Cherepova \cite{Ch99}   the author considered the heat volume potential acting on parabolic H\"older continuous functions that are allowed to blow up at the parabolic boundary. Finally,   Karazym and Suragan \cite{KaSu20} have considered the volume potential associated with a 
degenerate parabolic equation. However we note that, up to the author knowledge, no results for the heat volume potential on spaces of distributions are available in the literature.
 
 The results of the present paper continue the line of the works \cite{LaLu17, LaLu19, Lu20} on  the properties 
 of  integral operators of potential type appearing in the framework of parabolic theory.
Incidentally, we mention that our interest in proving these kind of mapping properties for the heat volume potential is also of technical nature since an  equation of  type \eqref{eq:heatf} arises when one tries to pull-back the 
heat equation to another domain  requiring only optimal regularity assumptions on the domains.
 In a subsequent paper by 
Dalla Riva and the author \cite{DaLu21} we will indeed use the results of the present paper to prove 
perturbation results for layer heat potentials.

The paper is organized as follows. In Section \ref{sec:sch} we introduce the functional spaces that we need, {\it i.e.}  parabolic Schauder spaces. Then, in Section \ref{sec:vp} we recall the definitions of the heat volume potential and of the Newtonian volume potential. Moreover, we show that the heat and Newtonian volume potentials coincide, up to the sign, whenever the density is time-independent. Section \ref{sec:main} contains the main result of the paper. More in detail, here we introduce the heat volume potential on  a  space of distributions that are the time derivative of H\"older continuous functions and we prove a mapping property of this operator. Finally, in Section 
\ref{sec:pb} we apply the results to the Dirichlet and Neumann problems. For the clarity of exposition, we have postponed a known result of functional analysis regarding quotient spaces to Appendix \ref{appA}

 
 \section{Schauder spaces} \label{sec:sch}
 Let $\Omega$ be a bounded open subset of $\mathbb{R}^n$. Let $k \in \mathbb{N}$ and 
 $\alpha \in \mathopen]0,1[$. For the definition of sets and functions of the Schauder class $C^{k,\alpha}$  we refer, \textit{e.g.}, to Gilbarg and
Trudinger~\cite{GiTr83}. Next we pass to recall the definitions of the parabolic analog of Schauder spaces.
Let 
 $T \in  \mathopen]-\infty,+\infty]$.   Let $\mathbb{D} \subseteq \mathbb{R}^n$. For the sake of brevity, we set 
\[
\mathbb{D}_T \equiv \overline{\mathopen ]-\infty,T\mathclose[} \times \mathbb{D}, 
\qquad \partial_T \mathbb{D} \equiv \overline{\mathopen ]-\infty,T\mathclose[} \times \partial\mathbb{D} .
\]
 We now introduce  the definition of an anisotropic H\"{o}lder space where  H\"{o}lder regularity with respect 
to time and  space directions can differ.
  Let $\alpha, \beta\in \mathopen]0,1[$.
 Then $C^{\alpha;\beta}(\mathbb{D}_T)$
denotes the space of bounded continuous functions $u$ from $\mathbb{D}_T$ to ${\mathbb{R}}$ such that
\begin{align*}
\|u\|_{C^{\alpha;\beta}(\mathbb{D}_T)}
\equiv &\sup_{ \mathbb{D}_T }|u|
+\sup_{\substack{t_{1},t_{2}\in  \overline{\mathopen]-\infty,T[}\\ t_{1}\neq t_{2}}   }\sup_{x \in \mathbb{D}}
\frac{|u(t_{1},x)  -u(t_{2},x)  | }{|t_{1}-t_{2}|^{\alpha}}\\ \nonumber
&+\sup_{t\in\overline{\mathopen]-\infty,T[}} 
\sup_{\substack{ x_1,x_2 \in \mathbb{D}\\ x_{1}\neq x_{2}}}
\frac{|u(t,x_1)  -u(t,x_2)  | }{|x_{1}-x_{2}|^{\beta}}<+\infty.
\end{align*}
We also denote by $C^{\alpha;0}(\mathbb{D}_T)$
the space of bounded continuous functions $u$ from $\mathbb{D}_T$ to ${\mathbb{R}}$ such that
\begin{align*}
\|u\|_{C^{\alpha;0}(\mathbb{D}_T)}
\equiv &\sup_{ \mathbb{D}_T }|u|
+\sup_{\substack{t_{1},t_{2}\in  \overline{\mathopen]-\infty,T[}\\ t_{1}\neq t_{2}}   }\sup_{x \in\mathbb{D}}
\frac{|u(t_{1},x)  -u(t_{2},x)  | }{|t_{1}-t_{2}|^{\alpha}}<+\infty.
\end{align*}

With the aim of considering boundary value problems, we will also need higher order parabolic H\"older space, {\it i.e.}  parabolic Schauder spaces. Let $\alpha \in \mathopen]0,1[$.  Let $\Omega$ be an open subset of $\mathbb{R}^n$.
We denote by $C^{\frac{1+\alpha}{2};1+\alpha}(\overline{\Omega_{T}})$ the space of bounded continuous functions $u$ from $\overline{{\Omega}_{T}}$ to ${\mathbb{R}}$ that are continuously differentiable with respect to the space variables and such that
\begin{align*}
\|u\|_{C^{\frac{1+\alpha}{2};1+\alpha}(\overline{\Omega_{T}})}
\equiv &\sup_{  \overline{{\Omega}_{T}} }|u|
+\sup_{\substack{t_{1},t_{2}\in  \overline{\mathopen]-\infty,T[}\\ t_{1}\neq t_{2}}   }\sup_{x \in \overline{\Omega}}
\frac{|u(t_{1},x)  -u(t_{2},x)  | }{|t_{1}-t_{2}|^{\frac{1+\alpha}{2}}}\\ \nonumber
&+\sum_{i=1}^n\|\partial_{x_i}u\|_{C^{\frac{\alpha}{2};\alpha}(\overline{\Omega_{T}})}<+\infty.
\end{align*}
If $\Omega$ is of class $C^{1,\alpha}$, by local parametrizations it is possible to naturally define the space $C^{\frac{1+\alpha}{2};1+\alpha}(\partial_T\Omega)$. We refer to Lady\v{z}enskaja,  Solonnikov, and Ural'ceva \cite{LaSoUr68}  and Lanza de Cristoforis and Luzzini \cite{LaLu17, LaLu19} for  more detailed definitions of
parabolic Schauder spaces.

Since we will consider the heat volume potential on a specific space of distributions, we need the following definition. We denote by $C^{-1+\alpha;\beta}(\overline{\Omega_{T}})$ the space of distributions 
in $\Omega_T$ that are the (distributional) time derivative of a function in $C^{\alpha;\beta}(\overline{\Omega_{T}})$, endowed with the quotient norm. That is 
\[
C^{-1+\alpha;\beta}(\overline{\Omega_{T}}) \equiv \left\{\partial_t u : u \in C^{\alpha;\beta}(\overline{\Omega_{T}}) \right\}
\]
and 
\[
\|f\|_{C^{-1+\alpha;\beta}(\overline{\Omega_{T}})} \equiv \inf\left\{\|u\|_{C^{\alpha;\beta}(\overline{\Omega_{T}})} : 
f = \partial_t u \right\}.
\]  
It can be easily seen that all the above spaces endowed with their respective norms are Banach spaces (also see Theorem \ref{thm:quot} for the case of $C^{-1+\alpha;\beta}(\overline{\Omega_{T}})$). Finally,  when $T >0$, with a subscript $0$ in the above spaces we mean the Banach subspace made of functions that are zero before zero. For example,
\[
C_0^{\alpha;\beta}(\mathbb{D}_T) \equiv \left\{u \in C^{\alpha;\beta}(\mathbb{D}_T) : 
u(t,\cdot)=0 \quad \forall t\leq 0 \right\}.
\]
The spaces $C^{\alpha;0}_0(\mathbb{D}_T)$, $C_0^{\frac{1+\alpha}{2};1+\alpha}(\overline{\Omega_{T}})$, and $C_0^{\frac{1+\alpha}{2};1+\alpha}(\partial_T\Omega)$
  can be defined in the same way. Similarly
\[
C_0^{-1+\alpha;\beta}(\overline{\Omega_{T}}) \equiv \left\{\partial_t u : u \in C^{\alpha;\beta}(\overline{\Omega_{T}}), \, \mathrm{supp}\,(\partial_t u) \subseteq [0,+\infty[ \right\}.
\] 
 \section{The heat volume potential}\label{sec:vp}
 Let $S_n: \mathbb{R}^{1+n} \setminus \{0,0\} \to \mathbb{R}$ denote the fundamental solution of the heat operator, that is 
\begin{equation*} 
S_{n}(t,x)\equiv
\begin{cases}
 \frac{1}{(4\pi t)^{\frac{n}{2}} }e^{-\frac{|x|^{2}}{4t}}&\quad{\mathrm{if}}\ (t,x) \in  \mathopen]0,+\infty\mathclose[  \times \mathbb{R}^n\,, 
 \\
 0 &\quad{\mathrm{if}}\  t \in  \left(\mathopen]-\infty,0] \times \mathbb{R}^n\right)\setminus \{(0,0)\}\,.
\end{cases}
\end{equation*}
As it is well know, $S_n \in C^\infty(\mathbb{R}^{1+n}\setminus \{(0,0)\})$ and solves the heat equation in 
$\mathbb{R}^{1+n}\setminus \{(0,0)\}$.
 We recall a known bound for $S_n$ which can be found {\it e.g.} in 
 Ladyzhenskaja, Solonnikov and Ural'ceva \cite[p. 274]{LaSoUr68}:
for all $\eta \in \mathbb{N}^n$ and for all $h\in\mathbb{N}$  there exists a constant
$C_{\eta,h}>0$ such that 
\begin{align}\label{boundS}
\left|D_x^\eta \partial_t^hS_n(t,x)\right| \leq C_{\eta,h}  t^{-\frac{n}{2}-\frac{|\eta|}{2}-h}e^{-\frac{|x|^2}{8t}}
\qquad \forall (t,x)\in   \mathopen]0,+\infty\mathclose[  \times \mathbb{R}^n.
\end{align}
Let $\Omega$ be a bounded open subset of $\mathbb{R}^n$ and $T \in \mathopen]-\infty,+\infty]$.  Let $x_0 \in \Omega$.  
If $f \in L^\infty( \Omega_T)$, we define the heat volume potential $P[f]$ to be the function from 
$\overline{\Omega_T}$ to 
$\mathbb{R}$ defined by 
\begin{equation}\label{def:volpot}
P[f](t,x) \equiv \int_{-\infty}^{+\infty}\int_\Omega \Big(S_n(t-\tau,x-y)- \delta_{2,n}S_n(-\tau,x_0-y)\Big)f(\tau,y)\,dyd\tau \qquad \forall (t,x) \in \overline{\Omega_T},
\end{equation}
where $\delta_{i,j}$ denotes the Kronecker delta.
 We note that the above definition, in the case $n=2$,  
depends on the choice of $x_0 \in \Omega$. Indeed, a different choice of $x_0$ would provide 
a volume potential that differs by a constant.
  However, if $T \in \mathopen ]0,+\infty]$ and 
$\mathrm{supp}\, f \subseteq \overline{[0,T[}\times \Omega$
(this is the case needed when one considers an initial-boundary value problem with initial condition at $t=0$),
 then the volume potential
$P[f]$ no longer depends on $x_0$ neither in the case $n=2$ and 
\begin{equation*}
P[f](t,x) = \int_{0}^t\int_\Omega S_n(t-\tau,x-y)f(\tau,y)\,dyd\tau \qquad \forall (t,x) \in \overline{\mathopen[0,T\mathclose[} \times \overline{\Omega},
\end{equation*}
 which is the classical definition of heat volume potential. The above definition with the term
 $\delta_{2,n}S_n(-\tau,x_0-y)$ is needed to avoid summability issues of the kernel as $\tau \to -\infty$ in the case $n=2$. Indeed $S_2(t-\tau,x-y)$ behaves as 
$(t-\tau)^{-1}$ as $\tau \to -\infty$, while $S_2(t-\tau,x-y)-  S_2(-\tau,x_0-y)$ does not have the same problem. 
To see this fact, we fix $\tau< \mathrm{min}\{0,t\}$ and $x,y \in \Omega$. One has 
 \begin{align*}
  \Big|S_2(t-\tau,x-y)&-  S_2(-\tau,x_0-y)\Big| \\
  \leq& \,\,\Big|S_2(t-\tau,x-y)-  S_2(t-\tau,x_0-y)\Big| \\
  &+\Big|S_2(t-\tau,x_0-y)-  S_2(-\tau,x_0-y)\Big|.
 \end{align*}
Then, if we denote by $\{e_j\}_{j=1,\ldots,n}$ the standard basis of $\mathbb{R}^n$,   the fundamental theorem of calculus and the estimates \eqref{boundS} for the fundamental solution $S_n$ imply that
  \begin{align*}
 \Big|S_2(t-\tau,x-y)&-  S_2(t-\tau,x_0-y)\Big| \\
  &\leq \sum_{j=1}^n|x_j-{x_0}_j|\int_0^1\Big|\partial_{x_j}S_2(t-\tau,\lambda x+(1-\lambda)x_0-y) \Big|\,d\lambda\\
   &\leq \sum_{j=1}^nC_{e_j,0}|x_j-{x_0}_j|\frac{1}{(t-\tau)^{\frac{3}{2}}}\int_0^1 e^{-\frac{|\lambda x+(1-\lambda)x_0-y|^2}{8(t-\tau)}}  \,d\lambda\\
   &\leq  \sum_{j=1}^nC_{e_j,0}|x_j-{x_0}_j|\frac{1}{(t-\tau)^{\frac{3}{2}}}
 \end{align*} 
 and that 
 \begin{align*}
 \Big|S_2(t-\tau,x_0-y)&-  S_2(-\tau,x_0-y)\Big| \\
 &\leq |t|\int_0^1\Big|\partial_t S_2(\lambda t-\tau,x_0-y)\Big|\,d\lambda\\
 &\leq C_{0,1}|t| \int_0^1 \frac{1}{|\lambda t-\tau|^{2}}e^{-\frac{|x_0-y|^2}{8(\lambda t -\tau)}}\,d\lambda\\
 &\leq C_{0,1}|t|  \frac{1}{|\min\{0,t\}-\tau|^{2}},
 \end{align*}
 which show that the kernel of \eqref{def:volpot} is summable for $\tau \to -\infty$.
 
We will also need  the Newtonian volume potential and then we recall that the fundamental solution of the Laplace equation is 
\begin{equation*}
\tilde S_n(x) \equiv 
\begin{cases}
\frac{1}{s_n}\log |x| &\forall x \in \mathbb{R}^n\setminus\{0\}, \quad  \mbox{ if } n=2,\\
\frac{1}{(2-n) s_n}|x|^{2-n}&\forall x \in \mathbb{R}^n\setminus\{0\},\quad \mbox{ if } n \geq 3,
\end{cases}
\end{equation*}
where 
\[
s_n \equiv \frac{2\pi^\frac{n}{2}}{\Gamma\left(\frac{n}{2}\right)}
\] denotes the $(n-1)$-dimensional measure of the unit sphere $\partial\mathbb{B}_n(0,1)$ and 
$\Gamma$ denotes the Euler Gamma function. If 
$h \in L^\infty(\Omega)$, the  harmonic volume potential $\tilde P[f]$ is the function from $\overline{\Omega}$ to 
$\mathbb{R}$ defined by 
\[
\tilde P[h](x) \equiv \int_\Omega \tilde S_n(x-y)h(y)\,dy \qquad \forall x \in \overline{\Omega}.
\]
Likewise other potential-type operators, the heat volume potential of an autonomous ({\it i.e.} time-independent) density   is autonomous and  coincides up to the sign with the corresponding Newtonian volume potential. That is, we have the following.
\begin{lemma}\label{lem:heatharm}
Let $\Omega$ be a bounded open subset of $\mathbb{R}^n$ and $T \in \mathopen]-\infty,+\infty]$.  Let $x_0 \in \Omega$.  Let 
$h \in L^\infty(\Omega)$. Then
\begin{equation*}
P[h](t,x) = -\big(\tilde P[h](x) - \delta_{2,n}\tilde P[h](x_0)\big)  \qquad \forall (t,x) \in \overline{\Omega_T}.
\end{equation*} 
\end{lemma}
\begin{proof}
We follow the lines of the proof of  \cite[Lemma A.3, Lemma A.4]{Lu19} where the analog relation between heat and harmonic layer potentials has been proved.
Let $(t,x)  \in \overline{\Omega_T}$.
We first consider the case $n=2$.  Then
\begin{align*}
P[h](t,x) =& \int_{-\infty}^{+\infty}\int_\Omega \Big(S_2(t-\tau,x-y)- S_2(-\tau,x_0-y)\Big)h(y)\,dyd\tau \\
=& \lim_{\sigma \to+ \infty}  \int_{-\sigma}^{+\infty} \int_{\Omega} \Big(S_2(t-\tau,x-y)- S_2(-\tau,x_0-y)\Big) h(y) \, dy d\tau.
\end{align*}
By the changes of variable $t-\tau = \frac{|x-y|^2}{4\xi}$ in the first term inside the integral and  $-\tau = \frac{|x_0-y|^2}{4\xi}$  in the  second term, we get 
\begin{align*}
P[h]&(t,x)\\
= & \lim_{\sigma \to +\infty}  \bigg\{ \int_{\Omega}\int_{\frac{|x-y|^2}{4(t+\sigma)}}^{+\infty} \frac{1}{4\pi \xi}e^{-\xi}h(y) \, d\xi dy -\int_{\Omega}\int_{\frac{|x_0-y|^2}{4\sigma}}^{+\infty} \frac{1}{4\pi \xi}e^{-\xi}h(y) \, d\xi dy \bigg\} \\
= & \lim_{\sigma \to +\infty}  \int_{\Omega}\int_{\frac{|x-y|^2}{4(t+\sigma)}}^{\frac{|x_0-y|^2}{4\sigma}}\frac{1}{4\pi \xi}e^{-\xi}h(y) \, d\xi dy.
\end{align*}
Let $g$ be the function from $\mathbb{R}$ to $\mathbb{R}$ defined by
\begin{equation*}
g(\xi) \equiv 
\begin{cases}
\frac{e^{-\xi}-1}{-\xi} &\mbox{ if } \xi \neq 0,\\
1   &\mbox{ if } \xi = 0.
\end{cases}
\end{equation*}
It is easy to see that $g$ is continuous in $\mathbb{R}$ and that 
\[
\xi^{-1} e^{-\xi} = \xi^{-1} - g(\xi) \qquad \forall \, \xi \in \mathbb{R} \setminus \{0\}
\] 
Accordingly, the dominated convergence theorem implies that
\begin{align*}
P[h](t,x) 
= & \lim_{\sigma \to +\infty}  \bigg\{ \int_{\Omega}\int_{\frac{|x-y|^2}{4(t+\sigma)}}^{\frac{|x_0-y|^2}{4\sigma}}\frac{1}{4\pi \xi} \, d\xi\,h(y)dy
-\int_{\Omega}\int_{\frac{|x-y|^2}{4(t+\sigma)}}^{\frac{|x_0-y|^2}{4\sigma}}\frac{g(\xi)}{4\pi} d\xi\,h(y)dy\bigg\}\\
= & \lim_{\sigma \to +\infty}  \bigg\{ \int_{\Omega}\frac{1}{4\pi }\log\left|\frac{|x_0-y|^2}{4\sigma}\frac{4(t+\sigma)}{|x-y|^2}\right|  h(y)\,dy-\int_{\Omega}\int_{\frac{|x-y|^2}{4(t+\sigma)}}^{\frac{|x_0-y|^2}{4\sigma}}\frac{g(\xi)}{4\pi} d\xi\,h(y)dy\bigg\}\\
= &  \int_{\Omega}\frac{1}{4\pi }\log\left(\frac{|x_0-y|^2}{|x-y|^2}\right)h(y)  \,dy\\
= &  \int_{\Omega}\frac{1}{2\pi }\log(|x_0-y|)  \, h(y)dy
- \int_{\Omega}\frac{1}{2\pi }\log(|x-y|)  \, h(y)dy\\
=&-\big(\tilde P[h](x) - \tilde P[h](x_0)\big).
\end{align*}
Next we pass to the case $n \geq 3$. By the change of variable 
$|x-y|^2s = 4(t-\tau)$ we have that 
\begin{align*}
P[h](t,x)   =&  \int_{-\infty}^{t} \int_{\Omega}\frac{1}{(4\pi (t-\tau))^\frac{n}{2}}e^{-\frac{|x-y|^2}{4(t-\tau)}}h(y)\, dy d\tau \\
=&\; \frac{1}{4\pi^\frac{n}{2}}\int_{0}^{+\infty}s^{-\frac{n}{2}}e^{-\frac{1}{s}}\,ds \int_{\Omega}\frac{1}{|x-y|^{n-2}}h(y)\, dy  \\
=&\;\frac{1}{4\pi^\frac{n}{2}}\Gamma \left(\frac{n}{2}-1\right) \int_{\Omega}\frac{1}{|x-y|^{n-2}}h(y)\, dy  \\
=&\;\frac{1}{(n-2) s_n}  \int_{\Omega}\frac{1}{|x-y|^{n-2}}h(y)\, dy \\
=&\;-\tilde P[h](x),
\end{align*}
which proves the statement.
\end{proof}

 
 \section{The heat volume potential on $C^{-1+\alpha,\beta}(\overline{\Omega_T})$}\label{sec:main}
In the present section we consider the action of $P[\cdot]$ in the space  $C^{-1+\alpha,\beta}(\overline{\Omega_T})$. Since 
 $f \in C^{-1+\alpha,\beta}(\overline{\Omega_T})$ is a distributions, and in particular is not a function  in 
$L^\infty(\Omega_T)$, we must specify what we mean by $P[f]$.
To this aim, we need  some preliminary results.
\begin{proposition}\label{prop:regB}
Let $\Omega$ be a bounded open subset of $\mathbb{R}^n$ and $T \in \mathopen]-\infty,+\infty]$. 
Let $\alpha \in \mathopen]0,1[$. Then the 
operator $B$ defined by 
\[
B[f](t,x) \equiv \int_{-\infty}^t\int_{\Omega} \partial_tS_n(t-\tau,x-y)(f(\tau,y)-f(t,y))\,dyd\tau \quad \forall (t,x)\in \overline{\Omega_T}
\]
is linear and continuous from $C^{\frac{1+\alpha}{2};0}(\overline{\Omega_T})$ to $C^{\frac{1+\alpha}{2};1+\alpha}(\overline{\Omega_T})$.
\end{proposition}
\begin{proof}
For convenience, for $\gamma >0$ we set 
\[
K_\gamma \equiv  \sup_{x \in \overline \Omega} \,\int_{\Omega}\frac{1}{|x-y|^{n-\gamma}}.
\]
It is easily seen that for r $\gamma >0$ one has
\begin{equation}\label{bk}
K_\gamma < +\infty.
\end{equation}
Let $(t,x) \in \overline{\Omega_T}$. By \eqref{boundS} there exists a constant $C_{0,1}>0$ such that
\begin{align*} 
 \int_{-\infty}^t \int_{\Omega} &\Big| \partial_t S_n(t-\tau,x-y)\big(f(\tau,y)-f(t,y)\big)\Big|\,dyd\tau \\\nonumber
&\leq C_{0,1}  \;\|f\|_{C^{\frac{1+\alpha}{2};0}(\overline{\Omega_T})} \int_{-\infty}^t \int_{\Omega} 
(t-\tau)^{-\frac{n}{2}-1+\frac{1+\alpha}{2}} e^{-\frac{|x-y|^2}{8(t-\tau)}}\;dyd\tau\\\nonumber
&= 8^{\frac{n-1-\alpha}{2}}C_{0,1}\|f\|_{ C^{\frac{1+\alpha}{2};0}(\overline{\Omega_T})} \int_{0}^{+\infty} s^{-\frac{n}{2}-1+\frac{1+\alpha}{2}} e^{-\frac{1}{s}}\;ds \int_{\Omega} \frac{1}{|x-y|^{n-1-\alpha}} \;dy\\\nonumber
&\leq 8^{\frac{n-1-\alpha}{2}}C_{0,1}K_{1+\alpha}\Gamma\left(\frac{n-1-\alpha}{2}\right) \;\|f\|_{ C^{\frac{1+\alpha}{2};0}(\overline{\Omega_T})}.
\end{align*}
Then, by the above inequality and by the Vitali convergence theorem,  $B[\cdot]$ is linear and continuous from $C^{\frac{1+\alpha}{2},\beta}(\overline{\Omega_T})$ to 
$C^{0}(\overline\Omega_T)$.
Next we take $t',t'' \in \mathopen ]-\infty,T[$, $t' < t''$,  $x \in \Omega$.  Then
\begin{align}\label{thm:mapprop1}
\big| B&[f](t',x) - B[f](t'',x) \big| \\ \nonumber
 \leq &\left| \int_{t'-2|t''-t'|}^{t'+2|t''-t'|} \int_{\Omega}\partial_t S_n(t'-\tau,x-y)\big(f(\tau,y)-f(t',y)\big)\,dyd\tau\right|\\ \nonumber
&+\left|  \int_{t''-3|t''-t'|}^{t''+3|t''-t'|} \int_{\Omega}\partial_t S_n(t''-\tau,x-y)\big(f(\tau,y)-f(t'',y)\big)\,dyd\tau \right|\\ \nonumber
&+\left|  \int_{-\infty}^{t'-2|t''-t'|} \int_{\Omega}\big(\partial_t S_n(t'-\tau,x-y)-\partial_t S_n(t''-\tau,x-y)\big)\big(f(\tau,y)-f(t',y)\big)\,dyd\tau \right|\\ \nonumber
&+\left|  \int_{-\infty}^{t'-2|t''-t'|} \int_{\Omega}\partial_t S_n(t''-\tau,x-y)\big(f(t'',y)-f(t',y)\big)\,dyd\tau \right|.
\end{align}
We begin considering the first term in the right hand side of \eqref{thm:mapprop1}. 
The bounds \eqref{boundS} on $S_n$  imply
  \begin{align*}
&\left| \int_{t'-2|t''-t'|}^{t' } \int_{\Omega}\partial_t S_n(t'-\tau,x-y)\big(f(\tau,y)-f(t',y)\big)\,dyd\tau\right|\\ \nonumber 
&\qquad \leq C_{0,1} \|f\|_{C^{\frac{1+\alpha}{2};0}(\overline{\Omega_T})} \int_{t'-2|t''-t'|}^{t' } \int_{\Omega}(t'-\tau)^{-\frac{n}{2}-\frac{1}{2}+\frac{\alpha}{2}}e^{-\frac{|x-y|^2}{8(t'-\tau)}}\,dyd\tau\\ \nonumber
&\qquad \leq C_{0,1} \|f\|_{C^{\frac{1+\alpha}{2};0}(\overline{\Omega_T})}\int_{t'-2|t''-t'|}^{t' }(t'-\tau)^{-\frac{n}{2}-\frac{1}{2}+\frac{\alpha}{2}}  \int_{\mathbb{R}^n}e^{-\frac{|x-y|^2}{8(t'-\tau)}}\,dyd\tau \\ \nonumber
&\qquad = (8\pi)^{\frac{n}{2}}C_{0,1}  \|f\|_{C^{\frac{1+\alpha}{2};0}(\overline{\Omega_T})}\int_{t'-2|t''-t'|}^{t' }(t'-\tau)^{-\frac{1}{2}+\frac{\alpha}{2}}  d\tau \\ \nonumber
&\qquad = (8\pi)^{\frac{n}{2}}\frac{2^{\frac{3+\alpha}{2}}}{1+\alpha} C_{0,1} \|f\|_{C^{\frac{1+\alpha}{2};0}(\overline{\Omega_T})}|t''-t'|^{\frac{1+\alpha}{2}}.
\end{align*}
The second  term in the right hand side of \eqref{thm:mapprop1} can be estimated in the same way. We then consider the third term. By \eqref{boundS} together with the mean value theorem one can see that there 
exists a constant $C'_{0,1}>0$ such that
\begin{equation}\label{in:ptS}
\left|\partial_t S_n(t'-\tau,x-y)-\partial_t S_n(t''-\tau,x-y) \right| \leq C'_{0,1}\frac{|t'-t''|}{|t'-\tau|^{\frac{n}{2}+2}}
e^{-\frac{|x-y|}{8(t'-\tau)}}
\end{equation}
for all $x,y \in \Omega$, $t'<t''$, $\tau< t'-2|t''-t'|$. For an explicit derivation of \eqref{in:ptS} we refer to 
Lanza de Cristoforis and Luzzini \cite[Lem. 4.3 (iii)]{LaLu17}. Then 
\begin{align*}
&\left|  \int_{-\infty}^{t'-2|t''-t'|} \int_{\Omega}\big(\partial_t S_n(t'-\tau,x-y)-\partial_t S_n(t''-\tau,x-y)\big)\big(f(\tau,y)-f(t',y)\big)\,dyd\tau \right|\\ \nonumber
&\qquad \leq C'_{0,1} \|f\|_{C^{\frac{1+\alpha}{2};0}(\overline{\Omega_T})}|t''-t'|
 \int_{-\infty}^{t'-2|t''-t'|} (t'-\tau)^{-\frac{n}{2}-2+\frac{1+\alpha}{2}}\int_{\Omega} e^{-\frac{|x-y|^2}{8(t'-\tau)}}\;dyd\tau \\ \nonumber
&\qquad \leq  C'_{0,1}\|f\|_{C^{\frac{1+\alpha}{2};0}(\overline{\Omega_T})}|t''-t'|
 \int_{-\infty}^{t'-2|t''-t'|} (t'-\tau)^{\frac{-n-3+\alpha}{2}}\int_{\mathbb{R}^n} e^{-\frac{|x-y|^2}{8(t'-\tau)}}\;dyd\tau \\ 
 &\qquad = (8\pi)^{\frac{n}{2}}  C'_{0,1}\|f\|_{C^{\frac{1+\alpha}{2};0}(\overline{\Omega_T})}|t''-t'|
 \int_{-\infty}^{t'-2|t''-t'|} (t'-\tau)^{\frac{-3+\alpha}{2}} d\tau \\ \nonumber
 &\qquad = \frac{2^{\frac{1+\alpha}{2}}}{1-\alpha} (8\pi)^{\frac{n}{2}} C'_{0,1}\|f\|_{C^{\frac{1+\alpha}{2};0}(\overline{\Omega_T})}
 |t''-t'|^{\frac{1+\alpha}{2}}.
\end{align*}
Next, we consider the last term in the right hand side of  \eqref{thm:mapprop1}.   
\begin{align*}
&\left|  \int_{-\infty}^{t'-2|t''-t'|} \int_{\Omega}\partial_t S_n(t''-\tau,x-y)\big(f(t'',y)-f(t',y)\big)\,dyd\tau \right|\\\nonumber
  &\qquad \leq \|f\|_{C^{\frac{1+\alpha}{2};0}(\overline{\Omega_T})}|t''-t'|^{\frac{1+\alpha}{2}}
   \int_{\Omega}\left|\int_{-\infty}^{t'-2|t''-t'|} \partial_t S_n(t''-\tau,x-y) \,d\tau\right|dy \\\nonumber
    &\qquad = \|f\|_{C^{\frac{1+\alpha}{2};0}(\overline{\Omega_T})}|t''-t'|^{\frac{1+\alpha}{2}}
    \int_{\Omega}\Big|S_n(3(t''-t'),x-y)\Big| \,dy \\\nonumber
     &\qquad = \|f\|_{C^{\frac{1+\alpha}{2};0}(\overline{\Omega_T})}|t''-t'|^{\frac{1+\alpha}{2}}
    \int_{\Omega}\frac{1}{(12 \pi (t''-t'))^\frac{n}{2}}e^{-\frac{|x-y|^2}{12(t''-t')}}\,dy \\\nonumber
     &\qquad \leq \|f\|_{C^{\frac{1+\alpha}{2};0}(\overline{\Omega_T})}|t''-t'|^{\frac{1+\alpha}{2}}.
\end{align*}
It remains to prove that for all $i \in \{1,\ldots,n\}$, the map $\partial_{x_i}B[\cdot]$ is linear and continuous from $C^{\frac{1+\alpha}{2},\beta}(\overline{\Omega_T})$ to 
$C^{\frac{\alpha}{2};\alpha}(\overline\Omega_T)$.  Let $(t,x) \in \Omega_T$. Let $i \in \{1,\ldots,n\}$.
 By classical differentiation theorems for  integrals depending on a parameter and by the bound  \eqref{boundS} there exists a constant 
$C_{e_i,1}>0$ such that   
\begin{align*}
|\partial_{x_i}B[f](t,x)|
\leq & \;C_{e_i,1}\|f\|_{C^{\frac{1+\alpha}{2};0}(\overline{\Omega_T})} \int_{-\infty}^t \int_{\Omega} 
(t-\tau)^{-\frac{n}{2}-\frac{3}{2}+\frac{1+\alpha}{2}} e^{-\frac{|x-y|^2}{8(t-\tau)}}\;dyd\tau\\
= & \;8^{\frac{n-1-\alpha}{2}}C_{e_i,1}\|f\|_{C^{\frac{1+\alpha}{2};0}(\overline{\Omega_T})} \int_{0}^{+\infty} u^{-\frac{n}{2}-1+\frac{\alpha}{2}} e^{-\frac{1}{u}}\;d\tau \int_{\Omega} \frac{1}{|x-y|^{n-\alpha}} \;dy\\
\leq & \;8^{\frac{n-1-\alpha}{2}}C_{e_i,1}K_\alpha\Gamma\left(\frac{n-\alpha}{2}\right)\|f\|_{C^{\frac{1+\alpha}{2};0}(\overline{\Omega_T})}.
\end{align*}
For what concerns the time $\frac{\alpha}{2}$-H\"older norm, it can be estimated exactly in the same way of the first part of the proof just by noting that the kernel is more singular by a term $(t-\tau)^{-\frac{1}{2}}$ but  we need to obtain an $\frac{\alpha}{2}$-H\"older regularity instead of an $\frac{1+\alpha}{2}$-regularity. 
Now we consider the spatial $\alpha$-H\"older regularity. Let $x',x'' \in \Omega$, $t \in \mathopen]-\infty,T[$.
By classical 
differentiation theorems for integrals depending on a parameter we have
\begin{align}\label{thm:mapprop2}
&|\partial_{x_i}B[f](t,x')-\partial_{x_i}B[f](t,x'')| \\\nonumber
& \quad = \bigg| \int_{0}^{+\infty} \int_{\Omega} \partial_{x_i}\partial_{t}S_n(\tau,x'-y)\big(f(t-\tau,y)-f(t,y)\big)\,dxd\tau  \\\nonumber
&\quad\,\,- \int_{0}^{+\infty} \int_{\Omega} \partial_{x_i}\partial_{t}S_n(\tau,x''-y)\big(f(t-\tau,y)-f(t,y)\big)\,dxd\tau  \bigg| \\\nonumber
& \quad \leq \left| \int_{0}^{|x'-x''|^2} \int_{\Omega} \partial_{x_i}\partial_{t}S_n(\tau,x'-y)\big(f(t-\tau,y)-f(t,y)\big)\,dxd\tau \right| \\\nonumber
& \quad \,\,+ \left| \int_{0}^{|x'-x''|^2} \int_{\Omega} \partial_{x_i}\partial_{t}S_n(\tau,x''-y)\big(f(t-\tau,y)-f(t,y)\big)\,dxd\tau \right| \\ \nonumber
& \quad\,\, + \left| \int_{|x'-x''|^2}^{+\infty} \int_{\Omega} \big(\partial_{x_i}\partial_{t}S_n(\tau,x'-y)-\partial_{x_i}\partial_{t}S_n(\tau,x''-y)\big)\big(f(t-\tau,y)-f(t,y)\big)\,dxd\tau \right|.
\end{align}
We consider the first term in the right hand side of \eqref{thm:mapprop2}. By the bound  \eqref{boundS} we have
\begin{align*}
& \left| \int_{0}^{|x'-x''|^2} \int_{\Omega} \partial_{x_i}\partial_{t}S_n(\tau,x'-y)\big(f(t-\tau,y)-f(t,y)\big)\,dyd\tau \right| \\
&\qquad \leq C_{e_i,1} \|f\|_{C^{\frac{1+\alpha}{2};0}(\overline{\Omega_T})} \int_{0}^{|x'-x''|^2} \int_{\Omega} \tau^{-\frac{n}{2}-\frac{3}{2}+\frac{1+\alpha}{2}} e^{-\frac{|x'-y|^2}{8\tau}}\,dyd\tau \\ 
&\qquad =(8\pi)^\frac{n}{2} C_{e_i,1} \|f\|_{C^{\frac{1+\alpha}{2};0}(\overline{\Omega_T})} \int_{0}^{|x'-x''|^2} \tau^{\frac{\alpha-2}{2}}  \,d\tau \\ 
&\qquad =\frac{2}{\alpha}(8\pi)^\frac{n}{2} C_{e_i,1} \|f\|_{C^{\frac{1+\alpha}{2};0}(\overline{\Omega_T})}  |x'-x''|^\alpha.
\end{align*}
The second term in the right hand side of \eqref{thm:mapprop2} con be estimated in the same way. Finally we consider the last term. The fundamental theorem of calculus and the bound  \eqref{boundS} imply
\begin{align*}
& \left| \int_{|x'-x''|^2}^{+\infty} \int_{\Omega} \big(\partial_{x_i}\partial_{t}S_n(\tau,x'-y)-\partial_{x_i}\partial_{t}S_n(\tau,x''-y)\big)\big(f(t-\tau,y)-f(t,y)\big)\,dyd\tau \right|\\
&\quad\leq  \|f\|_{C^{\frac{1+\alpha}{2};0}(\overline{\Omega_T})} \int_{|x'-x''|^2}^{+\infty} \tau^\frac{1+\alpha}{2} \int_{\Omega} \left|\big(\partial_{x_i}\partial_{t}S_n(\tau,x'-y)-\partial_{x_i}\partial_{t}S_n(\tau,x''-y)\right| \,dyd\tau\\
&\quad\leq  \|f\|_{C^{\frac{1+\alpha}{2};0}(\overline{\Omega_T})} \sum_{j=1}^n|x_j'-x''_j|\\
&\quad\qquad \times \int_0^1\int_{|x'-x''|^2}^{+\infty} \tau^\frac{1+\alpha}{2} \int_{\Omega} \left|\partial_{x_j}\partial_{x_i}\partial_{t}S_n(\tau,\lambda x'+ (1-\lambda) x''-y) \right| \,dyd\tau d\lambda\\
&\quad\leq  \|f\|_{C^{\frac{1+\alpha}{2};0}(\overline{\Omega_T})} \sum_{j=1}^nC_{e_i+e_j,1}|x_j'-x''_j|\\
&\quad\qquad \times \int_0^1\int_{|x'-x''|^2}^{+\infty} \tau^{\frac{n}{2}-2+\frac{1+\alpha}{2}}\int_{\Omega} 
e^{-\frac{|\lambda x'+ (1-\lambda) x''-y|^2}{8\tau}}  \,dyd\tau d\lambda\\
&\quad\leq (8\pi)^\frac{n}{2} \|f\|_{C^{\frac{1+\alpha}{2};0}(\overline{\Omega_T})}
\sum_{j=1}^nC_{e_i+e_j,1}|x_j'-x''_j|    \int_{|x'-x''|^2}^{+\infty}  \tau^{\frac{-3+\alpha}{2}}  \,d\tau \\
&\quad \leq \sum_{j=1}^nC_{e_i+e_j,1}\|f\|_{C^{\frac{1+\alpha}{2};0}(\overline{\Omega_T})}  |x'-x''| \int_{|x'-x''|^2}^{+\infty}  \tau^{\frac{\alpha-3}{2}} d\tau \\
&\quad= \frac{2}{1-\alpha} \sum_{j=1}^nC_{e_i+e_j,1}\|f\|_{C^{\frac{1+\alpha}{2};0}(\overline{\Omega_T})}  |x'-x''|^{\alpha}.
\end{align*}
\end{proof}
Next we show that the operator $B[f]$ of Proposition \ref{prop:regB} coincides with $\partial_t P[f]$
whenever $f \in  C^{\frac{1+\alpha}{2},\beta}(\overline{\Omega_T})$.
\begin{lemma}\label{lem:potrep}
Let $\Omega$ be a bounded open subset of $\mathbb{R}^n$ and $T \in \mathopen]-\infty,+\infty]$. Let $\alpha ,\beta \in \mathopen ]0,1[$. Let
$f \in C^{\frac{1+\alpha}{2},\beta}(\overline{\Omega_T})$. Then $P[f]$ is continuously differentiable with respect to $t$  and
\[
 \partial_tP[f](t,x)
= \int_{-\infty}^t\int_{\Omega} \partial_tS_n(t-\tau,x-y)(f(\tau,y)-f(t,y))\,dyd\tau \qquad \forall (t,x) \in {\Omega_T}.
\]
\end{lemma}
\begin{proof}
Let $(t,x) \in \Omega_T$. Since $f$ is 
$\beta$-H\"older continuous in space, by Friedman \cite[Thm. 9 p. 21]{Fr64} the volume potential $P[f]$ is continuously differentiable with respect to the time variable and two time 
continuously differentiable with respect to the space variables. Moreover
\begin{align*}
\partial _tP[f](t,x) = f(t,x) +
\Delta  P[f](t,x) .
\end{align*}
By the properties of the Newtonian volume potential (see {\it e.g.} Gilbarg and Trudinger 
\cite[\S 4.2]{GiTr83}), we have
\[
 \Delta \tilde P[f(t,\cdot)](x) = f(t,x)  .
\]
Since by Lemma \ref{lem:heatharm} the heat and Newtonian volume potential coincide up to the sign, we have that
\[
\partial _t P[f](t,x) = \Delta \Big(P[f](t,x) - P[f(t,\cdot)](t,x) \Big) = \Delta P[f-f(t,\cdot)](t,x).
\]
The bound \eqref{boundS} for the fundamental solution and \eqref{bk}  impliy that   
\begin{align*}
&\int_{-\infty}^t \int_{\Omega} \Big| \Delta S_n(t-\tau,x-y)\big(f(\tau,y)-f(t,y)\big)\Big|\,dyd\tau  \\\nonumber
&\quad\qquad\leq  \;\sum_{j=1}^nC_{2e_j,0}\|f\|_{C^{\frac{1+\alpha}{2};0}(\overline{\Omega_T})} \int_{-\infty}^t \int_{\Omega} 
(t-\tau)^{-\frac{n}{2}-1+\frac{1+\alpha}{2}} e^{-\frac{|x-y|^2}{8(t-\tau)}}\;dyd\tau\\\nonumber
&\quad\qquad\leq \;8^\frac{n-1+\alpha}{2} \sum_{j=1}^nC_{2e_j,0}\|f\|_{ C^{\frac{1+\alpha}{2};0}(\overline{\Omega_T})} \int_{0}^{+\infty} u^{-\frac{n}{2}-1+\frac{1+\alpha}{2}} e^{-\frac{1}{u}}\;d\tau \int_{\Omega} \frac{1}{|x-y|^{n-1-\alpha}} \;dy\\\nonumber
&\quad\qquad\leq \;8^\frac{n-1+\alpha}{2} \sum_{j=1}^nC_{2e_j,0}K_{1+\alpha}\Gamma\left(\frac{n-1-\alpha}{2}\right)\|f\|_{ C^{\frac{1+\alpha}{2};0}(\overline{\Omega_T})}.
\end{align*}
Accordingly the statement follows by standard differentiation theorems for integral depending on a parameter and by recalling that $S_n$ solves the heat equation in $\mathbb{R}^{1+n} \setminus \{(0,0)\}$.
\end{proof}

In order to define  the heat volume potential in $ C^{\frac{-1+\alpha}{2},\beta}(\overline{\Omega_T})$, which is a quotient space, we need to show that our definition is independent on the choice of the representative in the equivalence class.  To this aim, we need the following.
\begin{lemma}\label{lem:B0}
Let $\Omega$ be a bounded open subset of $\mathbb{R}^n$ and $T \in \mathopen]-\infty,+\infty]$. Let $\alpha ,\beta \in \mathopen ]0,1]$. Let
$f \in C^{\frac{1+\alpha}{2},\beta}(\overline{\Omega_T})$ be such that $\partial_t f =0$ in the sense of distributions. Then 
\begin{equation}\label{B0_1}
B[f](t,x)=\int_{-\infty}^t\int_{\Omega}\partial _tS_n(t-\tau,x-y)(f(\tau,y)-f(t,y))\,dyd\tau=0 \qquad \forall (t,x) \in 
\overline{\Omega_T}.
\end{equation}
\end{lemma}
\begin{proof}
Since $B[f]$ is continuous in $\overline{\Omega_T}$ by Proposition \ref{prop:regB}, it suffices to show equality 
\eqref{B0_1} in $\Omega_T$. Let $(t,x) \in \Omega_T$ be fixed. Since
$S_n(t-\tau,x-y)(f(\tau,y)-f(t,y))$ is continuous in $(\tau,y) \in \overline{\Omega_t} \setminus \{(t,x)\}$, it has a distributional $\tau$-derivative which, since $\partial_tf=0$ in the sense of distributions, equals
\[
g(\tau,y) \equiv -\partial _tS_n(t-\tau,x-y)(f(\tau,y)-f(t,y)) \qquad (\tau,y) \in {\Omega_t} \setminus \{(t,x)\}.
\]
Let  $\varepsilon>0$. Since the function 
$g(\tau,y)$ is continuous in $\overline{\Omega_{t-\varepsilon}}$, then
\begin{align*}
-\int_{-\infty}^{t-\varepsilon}\int_{\Omega} \partial_t&S_n(t-\tau,x-y)(f(\tau,y)-f(t,y))\,dyd\tau \\
&=\int_{\Omega} S_n(\varepsilon,x-y)(f(t-\varepsilon,y)-f(t,y))\,dy\\
&\quad \, -\lim_{\tau \to -\infty}\int_{\Omega} S_n(t-\tau,x-y)(f(\tau,y)-f(t,y))\,dy\\
&=\int_{\Omega} S_n(\varepsilon,x-y)(f(t-\varepsilon,y)-f(t,y))\,dy.
\end{align*}
For all $y \in \Omega$ we have
\begin{align*}
\Big| S_n(\varepsilon,x-y)&(f(t-\varepsilon,y)-f(t,y))\Big|\\
 &\leq \|f\|_{C^{\frac{1+\alpha}{2},\beta}(\overline{\Omega_T})}
\frac{1}{(4\pi)^{\frac{n}{2}}}\varepsilon^{\frac{-n+1+\alpha}{2}}e^{-\frac{|x-y|^2}{4\varepsilon}}\\
&\leq  \|f\|_{C^{\frac{1+\alpha}{2},\beta}(\overline{\Omega_T})}\left(\sup_{\xi >0}\xi^{\frac{-n+1+\alpha}{2}}e^{-\frac{1}{4\xi}}\right)\frac{1}{|x-y|^{n-1-\alpha}}.
\end{align*}
By the dominated convergence theorem we obtain $\eqref{B0_1}$ by letting $\varepsilon \to 0$.
\end{proof}
We are now ready to define the volume potential in $ C^{\frac{-1+\alpha}{2},\beta}(\overline{\Omega_T})$ as
\[
P[g](t,x) \equiv   \int_{-\infty}^t\int_{\Omega}\partial _tS_n(t-\tau,x-y)(f(\tau,y)-f(t,y))\,dyd\tau
\qquad \forall (t,x) \in \overline{\Omega_T},
\]
for all $g = \partial_t f \in C^{\frac{-1+\alpha}{2},\beta}(\overline{\Omega_T})$. Our main result is the following.
\begin{theorem}\label{thm:main}
Let $\Omega$ be a bounded open subset of $\mathbb{R}^n$ and $T \in \mathopen]-\infty,+\infty]$. Let $\alpha ,\beta \in \mathopen ]0,1[$. Then
\begin{itemize}
\item[i)] If $g \in C^{\frac{-1+\alpha}{2},\beta}(\overline{\Omega_T})$, then $\partial_t P[g] - \Delta P[g] = g$ 
in the sense of distributions in $\Omega_T$.
\item[ii)] $P[\cdot]$ is a bounded linear operator from  $C^{\frac{-1+\alpha}{2},\beta}(\overline{\Omega_T})$ to 
$C^{\frac{1+\alpha}{2},1+\alpha}(\overline{\Omega_T})$.
\end{itemize}
\end{theorem}
 \begin{proof}
 We start considering i). If $g \in C^{\frac{-1+\alpha}{2},\beta}(\overline{\Omega_T})$ and 
 $g = \partial_t f$ with $f \in  C^{\frac{1+\alpha}{2},\beta}(\overline{\Omega_T})$, then by Lemma
 \ref{lem:potrep} we have
 \[
 P[g]= \partial_tP[f] \qquad \forall (t,x) \in \overline{\Omega_T}.
 \]
 Hence
\begin{equation}\label{main1}
 \partial_t P[g] - \Delta P[g] = \partial_t \left( \partial_t P[f] - \Delta P[f]  \right)= \partial_t f = g,
 \end{equation}
 in the sense of distributions. The second equality in \eqref{main1} follows by classical results for the heat volume potential (see, {\it e.g.}, 
 Friedman \cite[Thm. 9, p.21]{Fr64}).
 
 Stament ii) simply follows by the definition of $P$, by Proposition \ref{prop:regB} and by Theorem  \ref{thm:quot} of the Appendix.
 \end{proof}


 \section{The Dirichlet and Neumann problems}\label{sec:pb}
 As an application, we show  the solvability of the Dirichlet and Neumann problem for  equation 
 \eqref{eq:heatf}. Let $\alpha, \beta \in \mathopen]0,1[$. Let $T>0$. Let $\Omega$ be a bounded open subset of $\mathbb{R}^n$ 
 of class $C^{1,\alpha}$. Let $f \in C_0^{\frac{1+\alpha}{2},\beta}(\overline{\Omega_T})$, 
 $g \in C_0^{\frac{1+\alpha}{2},1+\alpha}(\partial_T \Omega)$, and  $h \in C_0^{\frac{\alpha}{2},\alpha}( \partial_T \Omega)$. The Dirichlet problem for equation \eqref{eq:heatf} is
 \begin{equation}\label{prob:d}
 \begin{cases}
 \partial_t u- \Delta u = \partial_t f \qquad &\mbox{ in } ]0,T\mathclose[ \times \Omega,\\
 u = g \qquad &\mbox{ on } ]0,T\mathclose] \times \partial\Omega,\\
 u(0,\cdot) = 0 \qquad &\mbox{ in } \overline \Omega,
 \end{cases}
 \end{equation}
 while the Neumann problem reads 
 \begin{equation}\label{prob:n}
 \begin{cases}
 \partial_t v- \Delta v = \partial_t f \qquad &\mbox{ in } ]0,T\mathclose[ \times \Omega,\\
 \partial_\nu v = h \qquad &\mbox{ on } ]0,T\mathclose] \times \partial\Omega,\\
 v(0,\cdot) = 0 \qquad &\mbox{ in } \overline \Omega.
 \end{cases}
 \end{equation}
 We note that following the lines of the present section one can also show the solvability of problems for equation \eqref{eq:heatf} with boundary conditions other than Dirichlet and Neumann, as for example Robin or transmission boundary conditions. For the sake of simplicity, here we treat only problems  \eqref{prob:d} and 
 \eqref{prob:n}. We have the following existence and uniqueness result which is an immediate consequence of classical parabolic theory together with Theorem \ref{thm:main}.
 
 \begin{theorem}
 Let $\alpha, \beta \in \mathopen]0,1[$. Let $T >0$. Let $\Omega$ be a bounded open subset of $\mathbb{R}^n$ 
 of class $C^{1,\alpha}$. Let $f \in C_0^{\frac{1+\alpha}{2},\beta}(\overline{\Omega_T})$, 
 $g \in C_0^{\frac{1+\alpha}{2},1+\alpha}(\partial_T \Omega)$, and  $h \in C_0^{\frac{\alpha}{2},\alpha}( \partial_T\Omega)$. Then problem \eqref{prob:d} admits a unique solution 
 $u \in C_0^{\frac{1+\alpha}{2},1+\alpha}(\overline{\Omega}_T)$ and problem \eqref{prob:n} admits a unique solution 
 $v \in C_0^{\frac{1+\alpha}{2},1+\alpha}(\overline{\Omega}_T)$
 \end{theorem}
 
 \begin{proof}
 Since \eqref{prob:d} and \eqref{prob:n} are linear problems, the  uniqueness of their solutions in the space
 $C_0^{\frac{1+\alpha}{2},1+\alpha}( \overline{\Omega}_T)$ is well known (cf, {\it e.g.}, Friedman 
 \cite{Fr64}). 
 
 Next we pass to consider existence. By Theorem \ref{thm:main}, $P[\partial_t f] \in C_0^{\frac{1+\alpha}{2},1+\alpha}(\overline{\Omega}_T)$ and solves equation \eqref{eq:heatf}. Thus, the existence of a solution for problems  \eqref{prob:d} and \eqref{prob:n} can be reduced to the existence of a solution for 
 \begin{equation}\label{prob:da}
 \begin{cases}
 \partial_t \tilde u- \Delta \tilde u = 0\qquad &\mbox{ in } ]0,T\mathclose[ \times \Omega,\\
\tilde u = g-P[\partial_t f] \qquad &\mbox{ on } ]0,T\mathclose] \times \partial\Omega,\\
 \tilde u(0,\cdot) = 0 \qquad &\mbox{ in } \overline \Omega,
 \end{cases}
 \end{equation}
and for
 \begin{equation}\label{prob:na}
 \begin{cases}
 \partial_t \tilde v- \Delta\tilde  v = 0\qquad &\mbox{ in } ]0,T\mathclose[ \times \Omega,\\
 \partial_\nu \tilde v = h-\partial_\nu P[\partial_t f] \qquad &\mbox{ on } ]0,T\mathclose] \times \partial\Omega,\\
 \tilde v(0,\cdot) = 0 \qquad &\mbox{ in } \overline \Omega,
 \end{cases}
 \end{equation}
 respectively. It is classical that problems  \eqref{prob:da} and \eqref{prob:na} admit a solution 
 in the space $C_0^{\frac{1+\alpha}{2},1+\alpha}(  \overline{\Omega}_T)$. For a proof of this result based on potential theoretic methods we refer to Baderko \cite{Ba97} (see also Lunardi and Vespri \cite{LuVe91}).
 
 \end{proof}

 \begin{appendices}
\section{}\label{appA}
 In this Appendix we collect a well-know result in functional analysis regarding quotient spaces.
 For a proof we refer to Schaefer \cite[p. 42]{Sc71} (see also Dalla Riva, Lanza de Cristoforis and Musolino \cite{DaLaMu21}).
 \begin{theorem}\label{thm:quot}
 Let  $(X, \|\cdot\|_{X})$ be a normed space. Let $V$ be a closed subspace of $X$.
Then the norm on the quotient $X/V$ defined by
\[
\|[x]\|_{X/V} \equiv \inf_{v \in V}\|x+v\|_{X} \qquad \forall [x] \in X/V
\]
generates the quotient topology on $X/V$, {\it i.e} the strongest topology on $X/V$ such that the canonical projection $\pi$ of $X$ onto $X/V$ is continuous.

If $X$ is complete, then $(X/V, \|\cdot\|_{X/V})$ is complete. If $Y$ is a normed space and if $T$
 is a linear map from $X/V$ to $Y$, then $T$ is continuous if and only if $T \circ \pi$ is continuous.
 \end{theorem}
 \end{appendices}

\subsection*{Acknowledgment}
The author wish to thank Professor M. Lanza de Cristoforis and Professor M. Dalla Riva for many valuable 
comments and suggestions during the preparation of the paper.
The author is member of the `Gruppo Nazionale per l'Analisi Matematica, la Probabilit\`a e le loro Applicazioni' (GNAMPA) of the `Istituto Nazionale di Alta Matematica' (INdAM). The author also acknowledges the support of the Project BIRD191739/19 `Sensitivity analysis of partial differential equations in
the mathematical theory of electromagnetism' of the University of Padova.


\begin{thebibliography}{11}



\bibitem{Ba97}
E. A. Baderko,  Parabolic problems and boundary integral equations.  Math. Methods Appl. Sci. {\bf 20} (1997), no. 5, 449--459.

%
%
%
%

\bibitem{Ch99}
 M. F. Cherepova, On some properties of the parabolic potential of bulk masses. II. (Russian) Differ. Uravn. {\bf 36} (2000), no. 3, 408--414, 432; translation in Differ. Equ. {\bf 36} (2000), no. 3, 457--465.
%




\bibitem{DaLaMu15}
M. Dalla Riva, M. Lanza de Cristoforis and P. Musolino, Analytic dependence of volume potentials corresponding to parametric families of fundamental solutions. Integral Equations Operator Theory {\bf 82} (2015), no. 3, 371?393. 



\bibitem{DaLaMu20}
M. Dalla Riva, M. Lanza de Cristoforis and P. Musolino, Mapping properties of weakly singular periodic volume potentials in Roumieu classes.  
J. Integral Equations Appl. {\bf 32} (2020), no. 2, 129--149.

\bibitem{DaLaMu21}
M. Dalla Riva, M. Lanza de Cristoforis and P. Musolino.
Singularly Perturbed Boundary Value Problems: A Functional Analytic Approach. Springer International Publishing.

\bibitem{DaLu21}
M. Dalla Riva and   P. Luzzini,
Dependence of layer heat potentials upon perturbation of the support in the optimal H\"older setting. Preprint (2021).

%
%
%
%
%
%
 
  \bibitem{Fr64}
  A. Friedman,  Partial differential equations of parabolic type. Prentice-Hall, Inc., Englewood Cliffs, N.J. 1964 xiv+347 pp.

\bibitem{GiTr83}
D. Gilbarg , N.S. Trudinger,  Elliptic partial
differential equations of second order. Elliptic partial differential equations of second order. Reprint of the 1998 edition. Classics in Mathematics. Springer-Verlag, Berlin, 2001. xiv+517 pp. 


\bibitem{KaSu20}
M. Karazym, D. Suragan, Trace formulae of potentials for degenerate parabolic equations. Differential Integral Equations {\bf 33} (2020), no. 7-8, 337--360. 



%
%
 
 
  \bibitem{LaSoUr68}
O. A.~Lady\v{z}enskaja,  V.A.~Solonnikov, and N.N.~Ural'ceva,  Linear and quasilinear equations of parabolic type. (Russian) Translated from the Russian by S. Smith. Translations of Mathematical Monographs,   23 American Mathematical Society, Providence, R.I. 1968.

 
%
  
 
  \bibitem{LaLu17}
M. Lanza de Cristoforis and P. Luzzini, Time dependent boundary norms for kernels and regularizing properties of the double layer heat potential.  Eurasian Math. J. {\bf 8} (2017), no. 1, 76--118.
 
 
  \bibitem{LaLu19}
 M. Lanza de Cristoforis and P. Luzzini, Tangential derivatives and higher-order regularizing properties of the double layer heat potential. Analysis (Berlin) {\bf 38} (2018), no. 4, 167--193.


%

%
%

\bibitem{Li96}
G.M. Lieberman, Second order parabolic differential equations. World Scientific Publishing Co., Inc., River Edge, NJ, 1996. xii+439 pp.


\bibitem{Lu95}
A. Lunardi,  Analytic semigroups and optimal regularity in parabolic problems. [2013 reprint of the 1995 original] [MR1329547]. Modern Birkhäuser Classics. Birkhäuser/Springer Basel AG, Basel, 1995. xviii+424 pp. 

\bibitem{LuVe91}
A. Lunardi and V. Vespri, H\"older regularity in variational parabolic nonhomogeneous equations. J. Differential Equations {\bf 94} (1991), no. 1, 1--40. 


\bibitem{Lu19}
P. Luzzini, Regularizing properties of the double layer heat potential and shape analysis of a periodic problem. {\it Ph.D. Dissertation}, Universit\`a degli Studi di Padova (2019).


%

\bibitem{Lu20}
P. Luzzini,   Regularizing properties of space-periodic layer heat potentials and applications to boundary value problems in periodic domains. Math. Methods Appl. Sci. {\bf 43} (2020), no. 8, 5273--5294.

\bibitem{Mi65}
 C. Miranda, Sulle propriet\`a di regolarit\`a di certe trasformazioni integrali. (Italian) Atti Accad. Naz. Lincei Mem. Cl. Sci. Fis. Mat. Natur. Sez. Ia (8) 7 (1965), 303--336.

 
 
 



 
 \bibitem{Sc71}
  H. H. Schaefer,Topologica lvectors pace.Third printing corrected. GraduateTexts in 
  Mathematics, Vol. 3. Springer-Verlag, New York-Berlin, 1971.
 
%
%

\end{thebibliography}
\end{document}